\documentclass{amsart}
\usepackage{amsmath}
\usepackage{amssymb}
\usepackage[dvips]{graphicx}
\usepackage[latin1]{inputenc}
\usepackage{amsmath,xypic}
\usepackage{geometry}
\usepackage{amssymb,latexsym}

\newtheorem{thm}{Theora}[section]

\newtheorem{corollary}[thm]{Corollary}
\newtheorem{lemma}[thm]{Lemma}
\newtheorem{proposition}[thm]{Proposition}
\newtheorem{definition}[thm]{Definition}
\newtheorem{remark}[thm]{Remark}
\newtheorem{example}[thm]{Example}

\newtheorem{Definition and notations}[thm]{Definition and notations}

\newtheorem{Prop}[thm]{Proposition}

\newtheorem{Rem}[thm]{Remarque}

\title[A Survey on Non Characteristic Heisenberg Group Domains ]{A Survey on Non Characteristic Heisenberg Group Domains}

\author[Najoua Gamara- Nadia  Alamri]{}

\begin{document}
\maketitle

\centerline{ Najoua Gamara- Nadia  Alamri}
 \vspace*{0.5cm}
 {\bf Abstract:} In this work, we give a
survey on non characteristic domains of Heisenberg groups. We prove
that bounded domains which are diffeomorphic to the solid torus
having the center of the group as rotation axis, are non
characteristic. Then, we state the following conjecture : The
bounded non characteristic domains of the Heisenberg group of
dimension 1 are those diffeomorphic to a
solid torus having the center of the group as rotation axis.\\

\vspace*{0.2cm}

 {\bf Motivation:}  After the pioneered  works of D. Jerison and
J.M. Lee on the Yamabe problem on CR manifolds \cite{D.J}, and
precisely  during the last decade many works have been accomplished
on domains of the Heisenberg group and especially on non
characteristic ones. We can mention some of them: the resolution of
critical and  sub critical semi linear equations for the Kohn
Spencer Laplacian :
\begin{equation*}(P_{\varepsilon,f})
\left\{\begin{array}{llll} -\Delta_{H}u &=& K u^{3-\varepsilon}+f
&\text{ in }
\Omega\\ \;\;\;\;\;\;\;\;u &>&0 &\text{ in } \Omega\\
\;\;\;\;\;\;\;\;u &=&0 &\text{ on }
\partial\Omega\end{array}
\right.
\end{equation*}
 where $\Delta_{H}$ is the sublaplacian  on $\mathbb{H}^{1},$ $K$ is a
 $C^{3}$ positive
function defined on $\Omega,$ $f\in C(\overline{\Omega}),$ $f\geq0,$ $f\neq0$ and $\varepsilon\geq0.$\\
Different methods have been established to prove existence and
nonexistence results : variational methods
 \cite{SB,GC,IB.AC,NG.EL,NG.HG.AA}, sup and sub-solutions methods
 \cite{LB.MR.AGS,GL.JW,FU1},
 blow-up techniques  \cite{IB,IB.ICD,ICD.AC},  mean value formulas  \cite{EL.FU1,EL.FU2,FU2,FU3}, the moving
plane techniques  \cite{IB.JP}.  In  \cite{NG.HG.AA} to
investigate solutions of $(P_{\varepsilon,f})$ on  a bounded non
characteristic domain $\Omega$ of the Heisenberg group, the authors
use a method based on the study of the critical points at infinity
for the associated Euler-Lagrange functional and their effect on the
topology of its different level sets. Since in this case a crucial
role is played by the behavior of the Green function of the non characteristic domain  and its
regular part near the boundary, they  used the results of
\cite{NG.HG}, where a complete study of the Green Function, its
regular part  and their derivatives is given.\\  In \cite{NG.EL} for
Yamabe-type problems on Heisenberg group domains, $K=1$ and
$(P_{\varepsilon,f})=(P_{0,0})$,
 N. Garofallo and E. Lanconelli proved a nonexistence result for Yamabe-type problem ( the case $K=1$ and
$(P_{\varepsilon,f})=(P_{0,0})$), where  $\Omega$ is a
$H$-starshaped domain.
 In \cite{G.Citti.F.Uguzzoni},  G. Citti  and F. Uguzzoni  studied the CR version of
 a famous theorem due to A. Bahri and J.M. Coron \cite{A.B- J.M.C} and
 proved an existence result for Yamabe-type problem on Heisenberg
 group domains with nontrivial homology group.\\
In \cite{AM.VM},  the authors gave existence and multiplicity of
solutions for the cases $K=1, \varepsilon\geq0,\,\text{ in
}(P_{\varepsilon,f})$, on a bounded non characteristic domain using
the so called method of Liapunov-Schmidt reduction and variational
methods. The same authors with A.Pistoia in \cite{AM.VM.AP} proved
the existence of concentrating solutions for the slightly
sub-critical problem under a suitable assumption on $\partial\Omega$
and that the Robin's function of the domain has a non-degenerate
critical
point.\\
 Unfortunately few things are known
about the non characteristic sets of the Heisenberg group. For
instance, due to topological reasons every bounded C1 domain in the
Heisenberg group $\mathbb{H}^{n}$, whose boundary is homeomorphic to
the 2n-dimensional sphere $S^{2n}$, has non-empty characteristic
set. Recall that a basic result, due to Derridj \cite{De71},
\cite{De72}, shows that, at least from the measure theoretic point
of view, the set of characteristic points is not too big, more
precisely if $S$ is a $C^{1}$ surface of the Heisenberg group, the
standard surface measure of $S$ vanishes. We have to mention the
recent results of Balogh and Magnani, the first author \cite{Ba} has
proved that for a $C^{1}$ domain in the Heisenberg group
$\mathbb{H}^{n}$, the characteristic set has zero
$(Q-1)$-dimensional Haussdorf measure with respect to the
Carnot-Caratheodory distance of $\mathbb{H}^{n}$. Whereas Magnani
has extended Balogh's result to Carnot groups of step 2 in
\cite{Ma03(I)}, and further to groups of arbitrary step
\cite{Ma03(II)}.

 Typically, in the theory of
sub-elliptic equations, or in Carnot-Caratheodory Geometry,
characteristic points are present, take for example the more
simplest domains of $\mathbb{H}^{n}$:  Koranyi balls  have two
characteristic points. Consequently a domain with an empty
characteristic set must possess some special properties, either
geometric or topological. The goal of the present work is to provide
examples of domains such that their boundaries admit
  empty  characteristic sets.\\

\section{Introduction and Statement of main Results}
 In this work, we  will attempt to give a survey on non characteristic bounded domains of the Heisenberg group $\mathbb{H}^{1}$.
 As we will see, the construction of
bounded domains is  delicate, and  involves topology.

\par We begin by recalling that the  Heisenberg group $\mathbb{H}^{n},$ $(n\geq 1)$\;   is the Lie group whose underlying manifold is $\mathbb{C}^{n}\times\mathbb{R} =\mathbb{R}^{2n +1}$ and whose group law is given by
\begin{eqnarray*}
\tau_{(z',t')}(z,t) =(z',t')\cdot(z,t)= (x+x', y + y', t+t' + 2(<x,y'> - <x',y>))
\end{eqnarray*}
 where $<.,.>$ denotes the inner product in the euclidian space $\mathbb{R}^{n},$  $(z,t)=(x_{1},....,x_{n}, y_{1},....,y_{n},t)$ and
  $(z',t')=(x'_{1},....,x'_{n}, y'_{1},....,y'_{n},t')$.
The  Heisenberg group $\mathbb{H}^{n}$ is a Cauchy Riemann manifold
of hypersurface type. The horizontal distribution $\mathbf{H}$ of
$\mathbb{H}^{n},$ is spanned by the
following vector fields:\\
$X_{j}=\displaystyle\frac{\partial}{\partial
x_{j}}+2y_{j}\frac{\partial}{\partial t}\ \textrm{and}\
Y_{j}=\frac{\partial}{\partial y_{j}}-2x_{j}\frac{\partial}{\partial
t},\;\;j=1,..n$
and if we denote $\mathbf{T}=\displaystyle\frac{\partial}{\partial t},$ the tangent bundle of $\mathbb{H}^{n}$ has the following natural decomposition $$T\mathbb{H}^{n}=\mathbf{H}\bigoplus\mathbb{R}\mathbf{T}.$$  \\

Following the geometrical interpretation due to I. Piatetski-Shapiro, \cite{Ka}, ones can introduce the Heisenberg group
$\mathbb{H}^{n}$ using its identification with the boundary $M_{n}$ of the Siegel Domain:
 $$D_{n+1}=\{(\xi_{0}, \xi_{1},..., \xi_{n})= (\xi_{0}, \xi) \in \mathbb{C}\times\mathbb{C}^{n};  \sum_{1}^{n}|\xi_{j}|^{2}-Im \xi_{0}<0\}$$
$$M_{n}=\partial D_{n+1}=\{(\xi_{0}, \xi) \in \mathbb{C}\times\mathbb{C}^{n}; \sum_{1}^{n}|\xi_{j}|^{2}=Im \xi_{0}\}.$$
The Siegel domain $D_{n+1}$ is holomorphically equivalent to the
unit ball in $\mathbb{C}^{n+1}.$ The Heisenberg group
$\mathbb{H}^{n}$ acts on $\mathbb{C}^{n+1}$ by holomorphic affine
transformation which preserve $D_{n+1}$ and $M_{n}$ as follows: if
$(z,t)\in \mathbb{H}^{n}$ and $\xi\in \mathbb{C}^{n+1},$
$(z,t)\bullet \xi= \xi'$ where

\begin{eqnarray*}
\xi'_{0}
&=&\xi_{0}+t+i|z|^{2}+2i\sum_{1}^{n}\xi_{j} \bar{z}_{j}\\
\xi'_{j}&=&\xi_{j}+ z_{j},\;\;1\leq j\leq n.\\
\end{eqnarray*}
Since this action is transitive on $M_{n},$ the group
$\mathbb{H}^{n}$ is identified with $M_{n}$ via the correspondence:
$$(z,t)\leftrightarrow (z,t)\bullet 0= (t+i|z|^{2}, z_{1},......,
z_{n}).$$ Under this identification the CR structure on
$\mathbb{H}^{n}$ described above coincides with the CR structure on
$M_{n}$ induced from $\mathbb{C}^{n+1}.$\\

We will focus on the first component function of the correspondence
above which we denote  by $w=t+i|z|^{2}.$  In \cite{RH} R.Hladky has
stated that $w$ is a CR function and that  the domains of the
Heisenberg group $\mathbb{H}^{n}$ obtained as the product of
 a smoothly bounded precompact domain in the $2$-hyperbolic space and the unit sphere of dimension $2n-1$ for $n\geq2$, satisfy two important properties:\\
$(I)$ These domains have no characteristic boundary points.\\
$(II)$ These domains admit a smooth defining function  depending solely on the real and imaginary parts of the CR function $w.$\\

In this work, we will be interested on the study of the non
characteristic domains of  the $1$-dimensional Heisenberg group.

We consider  smooth bounded  domains of the upper half plane
$\mathbb{R}_{+}^{2},$.
$$ U\subset\mathbb{R}_{+}^{2}=\{(x,y) \in
\mathbb{R}^{2}  / \;\;\; y>0\}.$$ and their reciprocal images by $w$
in $\mathbb{H}^{1}:$
\begin{eqnarray}\label{nocenter}
 \Omega =\{(z,t) \in \mathbb{H}^{1}  / \;\;\;w=t+i|z|^{2} \in U\}=w^{-1}( U)
\end{eqnarray}
We prove  the following result:\\

$\bullet$ \textit{Let $\Omega$ be an open set of $\mathbb{H}^{1}$ and $U$
a smooth bounded  domain of the upper half plane as above, then  $\Omega$ is diffeomorphic to $U\times S^{1}.$}\\

We remark then that if the boundary of the domains $U$ intersect the
line $\{ y=0\},$ then the domains $\Omega$ satisfying
\eqref{nocenter} are characteristic sets of the Heisenberg group
$\mathbb{H}^{1}.$ Since in this work we are concerned with non
characteristic domains, we will take open sets $U$ of
$\mathbb{R}_{+}^{2}$ such that their boundaries
have no intersection with the line $\{ y=0\}.$ Under this  hypothesis, we prove the following results:\\

$\bullet$ \textit{Let $U$ be a convex domain of the upper half plane
as above and $A$ be any interior point of $U,$ then the closure of
$U$ is homeomorphic to a disc of center  $A.$}\\

$\bullet$ \textit{If $\Omega$ is a domain of $\mathbb{H}^{1}$
satisfying \eqref{nocenter} with $U$  convex  and $\partial
\Omega=w^{-1}(\partial U),$ then $\Omega$ is a non characteristic domain of the Heisenberg group $\mathbb{H}^{1}$.}\\

$\bullet$ \textit{The smoothly  bounded and open  domains  of
$\mathbb{H}^{1}$ which are diffeomorphic to the generalized solid
Torus of axis the center
of $\mathbb{H}^{1}$ are non characteristic.}\\

We conclude this survey by stating the following conjecture\\
$\bullet$ \textit{The only smoothly bounded  and non characteristic
domains of $\mathbb{H}^{1}$ are those  diffeomorphic to the
generalized solid Torus of revolution axis: the center of
$\mathbb{H}^{1}$.}

\section{Preliminaries}

The Heisenberg group $\mathbb{H}^{n}$ is the homogeneous Lie group
whose underlying manifold is $\mathbb{C}^{n}\times\mathbb{R}
=\mathbb{R}^{2n +1}$ and whose group law is given by
\begin{eqnarray*}
\tau_{(z,t)}(z',t') =(z,t)\cdot(z',t')= (x'+x, y' + y, t'+t + 2(<x',y> - <x,y'>))
\end{eqnarray*}
 where $<.,.>$ denotes the inner product in the euclidian space $\mathbb{R}^{n},$  $(z,t)=(x_{1},....,x_{n}, y_{1},....,y_{n},t)$ and $(z',t')=(x'_{1},....,x'_{n}, y'_{1},....,y'_{n},t')$.\\
The center of the Lie group $H^{n},$ is $Z=\{(0, t), 0\in
\mathbb{C}^{n}, t \in \mathbb{R} \}.$
 The $\mathbb{H}^{n}$-dilatations are the following transformations
\begin{equation*}
\delta_{\lambda}: \mathbb{H}^{n} \longrightarrow \mathbb{H}^{n}\;,\;
(z_{1},z_{2},....z_{n},t) \longmapsto(\lambda x, \lambda y,
\lambda^{2} t)\, , \,\,\text{ where }\lambda > 0.
\end{equation*}
The Jacobian determinant of $\delta_{\lambda}$ is $\lambda^{2n +2}$,
it yields that the homogeneous dimension of $\mathbb{H}^{n}$ is
 $2n+2$.\\ The homogeneous norm of the space is $ \rho(z_{1},z_{2},....z_{n},t)=
\big( (\sum_{j=1}^{n}(\mid z_{j} \mid^{4}  + t^{2}
\big)^{\frac{1}{4}} $ and the natural distance is accordingly
defined for
\begin{eqnarray*}
d((z,t);(z',t'))=
\rho((z,t)\cdot(z',t')^{-1}):=\rho((z,t)\cdot(-z',-t')).
\end{eqnarray*}
 The  Koranyi ball of center $\xi_{0}$ and radius $r$ is defined by  $B(\xi_{0},r)=\{\xi\in \mathbb{H}^n\; \rho(\xi_{0}\cdot\xi^{-1}) \leq
 r
\}.$ The vector fields
$$Z_{j}=\frac{\partial}{\partial z_{j}}+i \overline{z}_{j}\frac{\partial}{\partial t}, \;\;j=1,...n.$$
are left invariant with respect to the group law and homogenous of degree $-1$ with respect to the dilations.
 The space $\mathbf{T}_{1,0}=\text{span}\displaystyle\left\{ Z_{1},...........Z_{n}\right\}$ gives the CR structure of $\mathbb{H}^{n}$.\\
Let $\mathbf{H}=Re(\mathbf{T}_{1,0} + \mathbf{T}_{0,1})$ where $\mathbf{T}_{0,1}=\overline{\mathbf{T}_{1,0}}$. The  bundle $\mathbf{H}$
is called the horizontal distribution of $\mathbb{H}^{n}$ and it is spanned by
$X_{j}=\displaystyle\frac{\partial}{\partial x_{j}}+2y_{j}\frac{\partial}{\partial t}\ \textrm{and}
\  Y_{j}=\frac{\partial}{\partial y_{j}}-2x_{j}\frac{\partial}{\partial t},\;\;j=1,...n$.\\
The tangent bundle of $\mathbb{H}^{n}$ has the following natural decomposition $T\mathbb{H}^{n}=\mathbf{H}\bigoplus\mathbb{R}\mathbf{T}$,
 where $\mathbf{T}=\displaystyle\frac{\partial}{\partial t}.$\\
The complex structure on $\mathbb{H}^{n}$ is given by $$J:\mathbf{H}\longrightarrow \mathbf{H},\quad V+\overline{V}\longmapsto
 J(V+\overline{V})=i(V-\overline{V})\quad \forall\ V\in \mathbf{T}_{1,0}$$ and $(\mathbf{H},J)$ gives a real CR structure on $\mathbb{H}^{n}.$\\
The real 1-form $$\theta_{0}=dt+i \sum (z^{j}d\bar{z}_{j}- \bar{z}_{j}d z^{j})$$
annihilates $\mathbf{H},$ $ker(\theta_{0})=\mathbf{H},$ we take it to be the contact form for the CR structure of $\mathbb{H}^{n}$.\\
The subgradient or horizontal gradient of the Heisenberg group
$\mathbb{H}^n$ is given by $ \nabla_{\mathbb{H}^{n}}= (X_{1},....
X_{n},Y_{1},.... Y_{n}).$ The sublaplacian operator or the
Kohn-Spencer laplacian of $(\mathbb{H}^{n},\theta_{0})$  is given by
$$\Delta_{b}=-\sum_{j=1}^{n}\left(X_{j}^{2}+Y_{j}^{2}\right).$$

\begin{definition} Let $\Omega \subset \mathbb{H}^{n},$ be an open set.
 A $C^{m}$ defining function, $m \geq 1$, for $\Omega$
 is a real-valued $C^{m}$
function $\varphi$ defined on a neighborhood $U $ of the boundary of
$\Omega$ such that $\Omega \cap U =\{x \in U, \varphi(x) < 0 \}$,
$\partial\Omega \cap U =\{x \in U, \varphi(x) = 0 \}$
and $\nabla \varphi \neq 0$ on $\partial\Omega$, where $\nabla \varphi$ is the Euclidean gradient of $\varphi.$\\
 If $\Omega$ has a $C^{m}$ defining function, we say that it
 is a $C^{m}$ domain.
\end{definition}
\begin{definition}
Let $\Omega$ be a smooth domain of $\mathbb{H}^{n},$    $\xi_{0} \in
\partial\Omega,$ we say that a smooth function $\varphi$   describes the
boundary of $\Omega$ in a neighborhood of $\xi_{0}$ if there exist\\
$\begin{array}{rrcl}
\varphi :& B(\xi_{0})  & \longrightarrow& \;\mathbb{R};\\\\
\end{array}\;\;$
$\begin{array}{rrcl}
\varphi(\xi)&=0& \;\; \text{if}\;\; \xi \in B_{d}(\xi_{0}) \cap \partial\Omega\\
 \varphi(\xi)& <0 & \;\; \text{if}\;\; \xi \in B_{d}(\xi_{0}) \cap \Omega
\end{array}$\\
where  $B(\xi_{0})$ is a ball of $\mathbb{H}^{n}$ of center $\xi_{0}.$\\
\end{definition}

\begin{definition}
Let $\Omega$ be a smooth domain of $\mathbb{H}^{n},$    $\xi_{0} \in
\partial\Omega,$ we say that  $\xi_{0}$ is a
characteristic point of $\Omega,$ if
$\nabla_{\mathbb{H}^{n}}\varphi(\xi_{0})=0$,
 where  $\varphi$ is a smooth function which describes the boundary
of $\Omega$ in a neighborhood of $\xi_{0}.$
\end{definition}
We have also the following  equivalent definition
\begin{definition}
A point $\xi\in\partial\Omega$ is a characteristic point for
$\Omega$ if the boundary of $\Omega$ is tangent to the distribution
$H$ at $\xi$, i.e $T_{\xi}\partial\Omega=H_{\xi}.$
\end{definition}
\begin{remark}
Recall that the property  for a point of a given domain to be characteristic is independent of the choice of the function
 which locally describes this boundary since if $f_{1}$ and $f_{2}$
are two defining functions of the boundary of $\Omega$ near a point
$\xi_{0}$ then there exists a smooth function
 $h$ such that $f_{1}=h f_{2}$.\\
\end{remark}

At non-characteristic points the tangent space to $\partial\Omega$
intersects $H$ transversally with codimension $1.$ Thus,  at
characteristic points the horizontal normal to the boundary of the
domain could not be defined. This fact, makes analysis estimates
especially those related to  the resolution of critical and  sub
critical semi linear equations for the Kohn- Spencer laplacian
difficult near such points.

\begin{definition}
A Heisenberg group domain is said to be characteristic if its
boundary admits characteristic points, otherwise it is said to be
non characteristic.

\end{definition}
Recall that for a non characteristic point $\xi\in\partial\Omega,$
the intrinsic outer unit normal to $\partial\Omega$ at $\xi$ is
called the horizontal unit normal and it is given by
\begin{eqnarray}\label{nu}\overrightarrow{\nu}=\frac{\nabla_{\mathbb{H}^{n}}\varphi(\xi)}{|\nabla_{\mathbb{H}^{n}}\varphi(\xi)|_{H}}
\end{eqnarray}
\begin{remark}
For an interested reader detailed examples of characteristic and non
characteristic domains can be found in \cite{N.Alamri}.

\end{remark}
Let $U$ be a  regular bounded domain of the upper half plane   $\mathbb{R}_{+}^{2}$:\\
$$ U\subset\mathbb{R}_{+}^{2}=\{(x,y) \in \mathbb{R}^{2}  / \;\;\; y>0\}.$$
We consider the sub domains $\Omega$ of $\mathbb{H}^{1}$ satisfying
condition  \eqref{nocenter}. We have the following result.
\begin{Prop}\label{Diffeomorphic}
Let $\Omega$ be an open set of $\mathbb{H}^{1}$ with property
\eqref{nocenter}, then $\Omega$ is diffeomorphic to $U\times S^{1}.$
\end{Prop}
\textbf{Proof}: We  consider the following correspondence
 $$\begin{array}{rrclrrcl}
F :& \mathbb{H}^{1}\setminus Z & \longrightarrow& \;\; \mathbb{R}_{+}^{2} \times S^{1}\\
              & \;\;(z,t) & \longmapsto & \;\;(w(z,t), \frac{z}{|z|}).\\
\end{array}$$
   The tangent map of $F$ at any point $(z,t) \in \mathbb{H}^{1}\setminus Z$ is
$$\begin{array}{rrcl}
TF_{(z,t)} :& \mathbb{R}^{3} & \longrightarrow& \;\;\;\;\;\mathbb{R}^{2} \times TS^{1}_{\frac{z}{|z|}}\\
              & \;\;V=(v_{1},v_{2}, v_{3}) & \longmapsto & \;\;(v_{3}, [z\overline{(v_{1},v_{2})} +\overline{z}(v_{1},v_{2})]; \frac{(v_{1},v_{2})}{|z|} -\frac{1}{2}[z\overline{(v_{1},v_{2})} +\overline{z}(v_{1},v_{2})] \frac{z}{|z|^{3}}).
\end{array}$$
It is easy to check that $TF_{(z,t)}$ is an isomorphism, in fact it
is sufficient to prove that it is  one to one, which is the case
since if $V=(v_{1},v_{2}, v_{3})$ is in the kernel of  $TF_{(z,t)},$
we have
\begin{eqnarray*}
(v_{3}, [z\overline{(v_{1},v_{2})} +\overline{z}(v_{1},v_{2})];
\frac{(v_{1},v_{2})}{|z|} -\frac{1}{2}[z\overline{(v_{1},v_{2})}
+\overline{z}(v_{1},v_{2})] \frac{z}{|z|^{3}})&=&0
\end{eqnarray*}
it is straightforward that $v_{3}=0$ and after identification of the
second and third components of the left and the right hand side of
the above equation, we derive that $[z\overline{(v_{1},v_{2})}
+\overline{z}(v_{1},v_{2})]=0$ and $(v_{1},v_{2})=0,$ therefore
\begin{eqnarray*}
(v_{1},v_{2},v_{3})  &=&0.
\end{eqnarray*}
Which gives that  $F$ is a local diffeomorphism between  $\mathbb{H}^{1}\setminus Z$ and  $\mathbb{R}_{+}^{2}.$ \\
Furthermore $F$ is bijective, indeed\\
 - if  $(x+iy;Y) \in
\mathbb{R}_{+}^{2} \times S^{1},$ then $(x+iy;Y)=F(y^{\frac{1}{2}}Y;
 x),$ hence $F$ is onto.\\
 - $F$ is also one to one, since if
 $(t,|z|^{2},\frac{z}{|z|})=(t',|z'|^{2},\frac{z'}{|z'|}),$ then
 $t=t',$  $|z|=|z'|$ and $\frac{z}{|z|}=\frac{z'}{|z'|}$ which gives
 $0,$ $z$ and $z'$ are positively collinear with the same modulus so
 $z=z'.$\\
 So we have proved that $F$ is a diffeomorphism.
Since the open set $\Omega$ of $\mathbb{H}^{1}$ satisfies
\eqref{nocenter}, it is also an open set of $\mathbb{H}^{1}\setminus
Z.$ Now, if we consider the restriction of $F$ to $\Omega,$ we
deduce that $\Omega$ is diffeomorphic to its image  $F(\Omega),$
which is obviously equal to $U \times S^{1}.$ The proof is thereby
complete.


\begin{example}
We turn now to study a specific  example of domains of $\mathbb{R}_{+}^{2}.$ Let us consider the domain $U$ bounded by the geodesic
of $\mathbb{R}_{+}^{2}:$  the  half circle of center $O=(0,0)$ and radius $1$  based on the $x$-axis and the $x$-axis, more precisely
$U=\{(x,y),/; x^{2}+ y^{2}, \; y>0\}.$ By rotating $U$ around the  $x$-axis and  regarding the correspondence $F,$ we obtain a diffeomorphism
  between $U\times S^{1} $ and  the open unit Koranyi ball with the set $]-1, 1[\times \{0\}$  omitted. Hence, the open set $\Omega$ obtained
   has at least two characteristic points.
\begin{Rem}
The set $ \overline{U}\times S^{1}$ is diffeomorphic to the unit Koranyi ball.
\end{Rem}

\end{example}

What we can conclude here, is if the boundary of  $U$ intersects the
line $\{y=0\},$ the sets $\Omega$ satisfying \eqref{nocenter} admit
characteristic points on their boundaries. Since the aim of the
present work is to characterize the sets of $H^{1}$ without
characteristic boundary, we will limit our search to the set of  sub
domains $U$ of $\mathbb{R}_{+}^{2}$ such that their boundaries have
no intersection with the line $\{y=0\},$ which implies that the
boundaries of the sets $\Omega$ satisfying \eqref{nocenter} have no
intersection with the center $Z$ of $H^{1}.$

From now on, we will consider open sets $\Omega$ of $H^{1}$
satisfying property \eqref{nocenter} with  $ \partial U \cap
\{y=0\}= \emptyset.$ Furthermore, we suppose that
\begin{eqnarray}\label{boundaries}
\partial \Omega= w^{-1}(\partial U).
\end{eqnarray}
 Let $\Omega$ be  a bounded domain  of $\mathbb{H}^{1}$ satisfying
 \eqref{nocenter} and
\eqref{boundaries}.

We denote by $f$ the first coordinate function of $F$ and introduce
the following composition of maps:
$$\begin{array}{rrclrrcl}
\Psi :& H^{1}\setminus Z & \longrightarrow^{f} & \;\;\mathbb{R}_{+}^{2} &\longrightarrow^{\psi}&  \mathbb{R}_{+}^{2} \times \overline{\mathbb{R}_{+}^{2}}& \longrightarrow ^{\phi}& \;\;\mathbb{R}\\
              & \;\;(z,t) & \longmapsto & w(z,t)& \longmapsto   &   (w, \overline{w}) & \longmapsto & \phi (w, \overline{w})
\end{array}$$
where, $\phi$ is a smooth function and
$\overline{\mathbb{R}_{+}^{2}}$ denotes the conjugate set of
$\mathbb{R}_{+}^{2}.$\\
 We denote $\Psi_{\Omega}$, the restriction of $\Psi$ to $\Omega$:

$$\begin{array}{rrclrrcl}
\Psi_{\Omega} :& \Omega & \longrightarrow^{f} & \;\; U &\longrightarrow^{\psi}&  U \times \overline{U}& \longrightarrow ^{\phi}& \;\;\mathbb{R}\\
              & \;\;(z,t) & \longmapsto & w(z,t)& \longmapsto   &   (w, \overline{w}) & \longmapsto & \phi (w, \overline{w})
\end{array}$$
where, here $\overline{U}$ denotes the conjugate set of the open
$U.$
\begin{lemma}
Suppose that $\phi \circ \psi$ is a smooth  function which describes
the boundary of  $U$   in a neighborhood of $w_{0} \in  \partial U,$
such that $w_{0}=f(\xi_{0}),$  $\xi_{0}=(z_{0}, t_{0})\in
\partial \Omega.$ If $\Omega$ satisfies \eqref{nocenter} and \eqref{boundaries}, then $\Psi= \phi \circ \psi \circ f$ is a locally
defining function of the domain $\Omega,$ near $\xi_{0}.$
\end{lemma}

\begin{proposition}\label{tangent spaceof Omega}

Let  $U$ and  $\Omega$ be given  domains as above and $\xi_{0}\in
\partial \Omega, w_{0}=f(\xi_{0}).$  The tangent space $T_{\xi_{0}}\partial
\Omega$
 is characterized as follows:
\begin{eqnarray*}
T_{\xi_{0}}\partial \Omega&=&\{v=(v_{1}, v_{2},v_{3})
\in\mathbb{R}^{3} /(v_{3},x_{0}v_{1}+y_{0}v_{2})\;\; \in
T_{w_{0}}\partial U\}.
\end{eqnarray*}
\end{proposition}
\textbf{Proof}\\
The orthogonal space to $\nabla_{\xi_{0}}  \Psi $ is equal to the
tangent space $ T_{\xi_{0}}(\Omega(\xi_{0})\cap\partial
\Omega)=T_{\xi_{0}}\partial \Omega$  which is of dimension $2,$  it
is given by:
\begin{eqnarray*}
 T_{\xi_{0}}\partial \Omega&=&\{v\in\mathbb{R}^{3} / <\nabla_{\xi_{0}}  \Psi, v>=0\}\\
 &=&\{v\in\mathbb{R}^{3} /D f(\xi_{0})v\in Ker (D( \phi \circ \psi) (w_{0}))\}
\end{eqnarray*}
where $<, >$ denotes the Euclidian inner product and $Ker (D( \phi
\circ \psi) w_{0})$ is the kernel of the derivative of the function
$\phi \circ \psi$ at $w_{0}.$ Since  $\phi \circ \psi$ is a local
defining function of $\partial U$ at $w_{0},$  this kernel
 is equal to the tangent
space of $U(w_{0}) \cap
\partial U$ at $ w_{0}.$  So, it is equal to the
tangent space of $\partial U$ at $w_{0}.$
 Therefore, for  $v=(v_{1}, v_{2},v_{3}) \in \mathbb{R}^{3},$ the condition $D f(\xi_{0})v\in Ker (D( \phi \circ \psi) w_{0})$
  can be reformulated as follows: $(v_{3},x_{0}v_{1}+y_{0}v_{2})\; \in T_{w_{0}}\partial U,$ the result follows.\\



\section{Determination of the tangent space $T\partial \Omega_{\xi_{0}}$ }
In this section, we will  determine the tangent space
$T_{\xi_{0}}\partial \Omega$ for open bounded domains
$U\subset\mathbb{R}_{+}^{2}$ and  subdomains $\Omega$ of
$ \mathbb{H}^{1}\setminus Z$ satisfying conditions $\eqref{nocenter}$ and $\eqref{boundaries}.$\\

We suppose that $U$ is a convex set, hence the closure $K$ of $U$ is
a compact convex set of $\mathbb{R}_{+}^{2}$ hence of
$\mathbb{R}^{2}.$ Let $A=(a_{1},a_{2})$ be an interior point of $U,$
we have the following result
\begin{proposition}\label{convex}
The compact convex set $K$ is homeomorphic to a disc of center  $A.$
\end{proposition}
\textbf{Proof:}\\
Since $U$ is an open set of $\mathbb{R}^{2},$ it contains  an open disc of center $A$ and radius $r,$  denote it by $D_{O}(A,r).$
Therefore the closed disc $D(A,r)\subset \bar{U}.$
Next, we will consider the following map which sends the boundary of $U$ on the boundary of the disc $D(A,r)$
$$\begin{array}{rrcl}
g :& \partial U  & \longrightarrow& \;S^{1}(A,r) \\
   &X& \;\mapsto &\;\; A + r\displaystyle\frac{X-A}{\|X-A\|}
\end{array}$$\\
 $g$ is one to one  since if $X$ and $X'$ are points of the boundary such that $g(X)=g(X'), $
  then $\displaystyle\frac{X-A}{\|X-A\|}=\displaystyle\frac{X'-A}{\|X'-A\|}, $ suppose that
   $X\neq X'$ which implies that  the points $X, X'$ and $A$ are collinear. Suppose for example
    that $X'$ is between $A$ and $X$ and consider the homothetic transformation  $H$ of $\mathbb{R}^{2}$
    of center $X'$ and ratio $h$ which sends $A$ to $X.$ Let $X_{0}\in D(A,r)$ and $Y_{0}=H(X_{0}),$ since
     $K$ is convex, the segment $[X_{0}, X']\subset K$ but $Y_{0}\in [X_{0}, X'],$ therefore $Y_{0}\in K.$
     Hence, the image $H(D(A,r))$ is a disc of center $X$ included in $K,$ it yields that $X$ is an interior
     point of $K$ but $X\in\partial U=\partial K$ which is a contradiction.
The map $g$ is onto:\\
Let $M$ be a boundary point of $K,$ since $K$ is convex and  $A\in
U\subset K$ the segment $[A, M]\in K.$ Let $\omega\in S^{1}(A,r),$
the intersection $\{A+ t(\omega - A), t\geq0\} \cap K$ of the half
line of $\mathbb{R}^{2}$  of origin $A$ containing $\omega$ with
$K,$  is compact and convex in $\{A+ t(\omega - A), t\geq0\},$ hence
it is equal to a segment $[A, A+ a(\omega)(\omega -A)]$ where
$a(\omega)>0$ is a constant depending on $\omega.$ The point $A+
a(\omega) (\omega -A)\in \partial K,$ indeed the sequence of points
$X_{m}=(A+ (a(\omega)+\frac{1}{m}) (\omega-A)),$ $m\geq1$ is in the
complement of $K$ in $\mathbb{R}^{2}$ and converges to $A
+a(\omega)(\omega -A)$ and $g(A+ a(\omega)( \omega-A))=\omega.$
Since the function $g$ is continuous and closed, we deduce that $g$
is an homeomorphism between $\partial K$ and $S^{1}(A,r).$ Now
define:

$$\begin{array}{rrclrrcl}
\mathfrak{G} :& D(A,r) & \longrightarrow & \;\; K \\
              & \;\;Y & \longmapsto &\Big\{\begin{array}{c}
                                                        \qquad \qquad  A \qquad\qquad\qquad\qquad\qquad \qquad \qquad \quad \text{if}\;\; Y= A\\
                                                       A+ \displaystyle\frac{\|A-Y\|}{\|A-Pr_{S^{1}}(Y)\|} [g^{-1}(Pr_{S^{1}}(Y))-A]\qquad \;\;\text{if}\;\;Y\neq A
                                                      \end{array}\end{array}$$
We claim that  $\mathfrak{G}$ is a homeomorphism between $D(A,r)$
and $K.$ It is obvious that the function $\mathfrak{G}$ is
continuous for $Y\neq A$ and the continuity in $A$ follows from the
compactness of $\partial K.$ On the other hand, $\mathfrak{G}$ is
bijective with inverse map

$$\begin{array}{rrclrrcl}
\mathfrak{H} :& K & \longrightarrow & \;\; D(A,r) \\
              & \;\;X & \longmapsto &\Big\{\begin{array}{c}
                                                        \qquad \qquad  A \qquad\qquad\qquad\qquad\qquad \qquad \qquad \quad \text{if}\;\; X= A\\
                                                       A+ \displaystyle\frac{\|A-X\|}{\|A-Pr_{\partial K}(X)\|} [g(Pr_{\partial K}\;X))-
                                                       A]\qquad \;\;\text{if}\;\;X\neq A
\end{array}
\end{array}$$
where $Pr_{\partial K} $ and $Pr_{S^{1}},$ are respectively the
projection map on the boundary of $K$ and on the boundary of the
circle $S^{1}(A,r).$\\
 To complete the proof, it remains
to show that $\mathfrak{G}$ is closed, this fact follows using the
compactness of the disc $D(A,r).$

\begin{proposition}\label{Non Characteristic}
Let $U$ be an open convex and bounded subset of the half upper plane
such that the boundary of $U$ does not intersect the line $\{y=0\}.$
If $\Omega$ is the domain of $\mathbb{H}^{1}$ satisfying
\eqref{nocenter} and \eqref{boundaries}, then
 $\Omega$ is a non characteristic domain
of $\mathbb{H}^{1}.$
\end{proposition}
\textbf{Proof}\\
 Using  proposition \ref{convex}, $U$
is diffeomorphic to the open disc $D_{o}(A,r)=\{M \in\mathbb{R}^{2}/
\|AM\|< r\}.$ Suppose that $A=(a_{1}, a_{2}),$ then $S^{1}(A,r)=
\zeta^{-1}(0)$ where $\zeta$ is the function

$$\begin{array}{rrclrrcl}
\zeta :& \mathbb{R}^{2} & \longrightarrow & \;\mathbb{R} \\
              & \;(X,Y) & \longmapsto & (X-a_{1})^{2}+(Y-a_{2})^{2} -r^{2}
                                                      \end{array}$$

Therefore, the tangent space to the circle $S^{1}(A,r)$ at the point
$(X,Y),$ $T_{(X,Y)}S^{1}(A,r)$ is equal to the vector space $<
\overrightarrow{AM}>^{\perp}.$ Since for a given  $\xi_{0}\in
\partial \Omega$ and  $w_{0}=f(\xi_{0}),$  we have
\begin{eqnarray*}
T_{\xi_{0}}\partial \Omega&=&\{v=(v_{1}, v_{2},v_{3}) \in\mathbb{R}^{3} /(v_{3},x_{0}v_{1}+y_{0}v_{2})\;\; \in T_{w_{0}}\partial U\}\\
\end{eqnarray*}
 we obtain\\
$T_{\xi_{0}}\partial \Omega=\{v=(v_{1}, v_{2},v_{3}) \in\mathbb{R}^{3} / v_{3}(t_{0}-a_{1})+ 2 (x_{0}v_{1}+y_{0}v_{2}) (x_{0}^{2}+ y_{0}^{2} -a_{2})= 0 \}.$\\

So, $v=(v_{1}, v_{2},v_{3})$ is in the tangent space of $\partial
\Omega$ at $\xi_{0}=(z_{0}, t_{0})$ if and only if the following
equation is satisfied
\begin{eqnarray}\label{tangency} v_{3}(t_{0}-a_{1})+ 2 (x_{0}v_{1}+y_{0}v_{2}) (x_{0}^{2}+ y_{0}^{2}
-a_{2})= 0.
\end{eqnarray}
 If $t_{0}=a_1$ then $(x_0^2+y_0^2) \neq a_2,$
otherwise $f(x_0,y_0,t_0)= (a_1,a_2)$ which is  an  interior point
of $U.$ For the same reason if $x_0^2+y_0^2=a_2,$ then $t_0 \neq
a_1,$ .\\
  We denote the horizontal tangent space to $\partial \Omega$ at $\xi_0$ by $\mathfrak{h}_{\xi_{0}},
  $ and
  let $V=(v_1,v_2,v_3)\in T_{\xi_{0}}\partial \Omega \cap \mathfrak{h}_{\xi_{0}}.$  $(X_{\xi_{0}}, Y_{\xi_{0}})$  is a basis of
   $\mathfrak{h}_{\xi_{0}}$: $$X_{\xi_{0}}(\xi)= \frac{\partial}{\partial x}+ 2y_0 \frac{\partial}{\partial t}$$
  $$Y_{\xi_{0}}(\xi)= \frac{\partial}{\partial y}- 2x_0 \frac{\partial}{\partial t}$$
   for $\xi=(x,y,t)\in \mathbb{H}^{1}.$
We have    $$V = v_1 \frac{\partial}{\partial x|_{\xi_{0}}}+v_2
\frac{\partial}{\partial y|_{\xi_{0}}} +v_3 \frac{\partial}{\partial
t|_{\xi_{0}}}$$
  $$ = \alpha\, X_{\xi_{0}}+ \beta \, Y_{\xi_{0}} \, \alpha, \beta \in \mathbb{R} .$$ It yields
   $$V= \alpha \frac{\partial}{\partial x|_{\xi_{0}}} + 2\alpha y_0 \frac{\partial}{\partial t|_{\xi_{0}}} +
   \beta \frac{\partial}{\partial y|_{\xi_{0}}}-2 \beta x_0 \frac{\partial}{\partial t|_{\xi_{0}}}.$$
So $$V= \alpha \frac{\partial}{\partial x|_{\xi_{0}}}+ \beta
\frac{\partial}{\partial y|_{\xi_{0}}} + 2(\alpha y_0 - \beta
x_0)\frac{\partial}{\partial t|_{\xi_{0}}}.$$
 Therefore, if  $ V=(v_1,v_2,v_3) \in
T_{\xi_{0}}\partial \Omega \cap \mathfrak{h}_{\xi_{0}},$ we obtain\\
$$ v_1 = \alpha,\,\, v_2 = \beta $$
and
\begin{eqnarray}\label{intersection}
v_3 &=2(\alpha y_0 - \beta x_0)= 2(v_1 y_0 - v_2 x_0)&
\end{eqnarray}
Three cases may occur:\begin{enumerate}
\item \textbf{First case:}\\
 If $t_{0}=a_1$ from equation \eqref{tangency}, we deduce that $v_1x_0+v_2y_0= 0,$ since $(x_0^2+y_0^2) \neq a_2,$
  hence $(v_1,v_2)$ is orthogonal to
  $(x_0,y_0),$ then $(v_1,v_2)$ is of the form: $(v_1,v_2)= a (-y_0 ,x_0),$ where
$a$ is a real constant.  Thus, we deduce the final form for the
vector $V=(v_1,v_2,v_3)$ in $T_{\xi_{0}}\partial \Omega \cap
\mathfrak{h}_{\xi_{0}}$ :$$V =a(-y_0,x_0,-2(v_0^2+y_0^2)).$$
 It yields that\\ $ T_{\xi_{0}}\partial \Omega \cap \mathfrak{h}_{\xi_{0}} =\{ a(-y_0, x_0,-2(x_0^2+y_0^2)), a \in \mathbb{R} \}$ which
  is a vector space of dimension $1.$ Therefore,
$\xi_{0}$ is not a characteristic point for $\Omega.$
 \item \textbf{Second case:}\\
 If $x_0^2+y_0^2=a_2,$ then from equation \eqref{tangency}, we deduce that $v_3 (t_0-a_1)=0$ since $t_0\neq a_1 then  v_3=0.$
 In this case, if $ V=(v_1,v_2,v_3) \in T_{\xi_{0}}\partial \Omega \cap \mathfrak{h}_{\xi_{0}},$  $V=(v_1,v_2,0).$
 But, we have from \eqref{intersection}, $v_3=2(v_1 y_0 - v_2 x_0),$ so $v_1\,y_0=v_2\,x_0.$ \\
 If $y_0\neq 0$ then, $v_1=v_2\displaystyle\frac{x_0}{y_0}$ and $V=(v_2\displaystyle\frac{x_0}{y_0}, v_2,\,0).$
 It yields that\\
   $T_{\xi_{0}}\partial \Omega \cap \mathfrak{h}_{\xi_{0}}=
  \{ v_2(\displaystyle\frac{x_0}{y_0},\,1,\,0), v_2 \in \mathbb{R} \},$
 which is a vector space of dimension $1.$\\
 If $x_0\neq 0,$ then  $v_2=v_1\displaystyle\frac{y_0}{x_0}$ and $V=(v_1, v_1\displaystyle\frac{y_0}{x_0},\,0).$ In this
 case \\
 $T_{\xi_{0}}\partial \Omega \cap \mathfrak{h}_{\xi_{0}}= \{ v_1(1,\,\displaystyle\frac{y_0}{x_0},\,0), v_1 \in \mathbb{R} \},$ which
 is also  a vector space of dimension $1.$\\
The case $y_0=x_0=0$ could not occur since by hypothesis
$$f(x_0,y_0,t_0)= w_0 = t_0+i|z_0|^2 \in \partial U$$ and
$$\partial U \cap \{ y=0 \} =\emptyset$$
so $|z_{0}|^2 \neq 0.$\\
\item \textbf{Third case}\\
If $t_0\neq a_1 $ and $ x_0^2+y_0^2\neq a_2$, we obtain  for
$V=(v_1,v_2,v_2)\in T_{\xi_{0}}\partial \Omega \cap
\mathfrak{h}_{\xi_{0}}$ using \eqref{tangency} and
\eqref{intersection}, the following equations

$$2x_0(x_0^2+y_0^2 -a_2)v_1 +2y_0(x_0^2+y_0^2 -a_2)v_2 + (t_0-a_1)v_3=0$$
$$2y_0v_1-2x_0v_2- v_3=0.$$
If we denote $x_0^2+y_0^2 -a_2$ by $B$ and $t_0-a_1$ by $C,$ we have

$$2x_0 B v_1 +2y_0 B v_2 + C (2y_0v_1-2x_0v_2)=0$$ and
$$2y_0v_1-2x_0v_2- v_3=0.$$
It yields that
$$[2x_0 B+2y_0 C]v_1+[2y_0 B-2x_0 C]v_2=0$$
and $$2y_0v_1-2x_0v_2- v_3=0.$$ This system of two equations is the
intersection of two  planes:
$$P_1: [B\,x_0+C\,y_0]v_1+[B\,y_0-C\,x_0]v_2 =0.$$
$$P_2: 2y_0\,v_1-2x_0\,v_2- v_3=0.$$
The planes $P_1$ and $P_2$ have respectively the following normal
vectors:$$N_1=([B\,x_0  +C\,y_0 ],[B\,y_0  -C\,x_0 ],0)$$ and
$$N_2= (2y_0,-2x_0 ,-1).$$ Since, $N_1$ is not parallel to $ N_2,$
we deduce that the  intersection of $P_1$ and $P_2$ is a vector
space of dimension $1.$ 
The result follows.
\end{enumerate}

\begin{remark}
The set $U\times S^{1}$ has a natural structure of a smooth manifold
homeomorphic to the open solid torus $D(A,r)\times S^{1}.$ It is
known that in dimension $1, 2 $ and $3$ any pair of homeomorphic
smooth manifolds are diffeomorphic.
\end{remark}

As a consequence of Proposition \ref{convex} and proposition
\ref{Non Characteristic} we have  the following result
\begin{corollary}
The open smoothly bounded domains of the Heisenberg group
$\mathbb{H}^{1}$ which are diffeomorphic to the solid generalized
open torus having as revolution axis the center of the group are non
characteristic.
\end{corollary}
Next, we give an example of a non characteristic domain of the Heisenberg group $\mathbb{H}^{1}.$
\begin{example}
We consider: $\Omega=\{(z,t)\in \mathbb{H}^{1}, |w-(2i+1)|< 1\}$.
The  set $\Omega$ is the open  Heisenberg solid torus $D_{o}\times
S^{1},$ having the center of $\mathbb{H}^{1}$ as revolution axis
where $D_{o}$ is the open disc of $\mathbb{R}^{2}_{+}$ of center $(1,2)$
and radius $1.$ 
\end{example}
And we conclude this survey by stating the following\\
\textbf{ Conjecture} \\
The only smoothly bounded and  non characteristic domains of the
Heisenberg group  $\mathbb{H}^{1}$ are those  diffeomorphic to the
generalized solid torus of revolution axis, the center of the group.


\end{document}